\documentclass[twoside, 11pt]{article}
\usepackage{mathrsfs}
\usepackage{amssymb, amsmath, mathrsfs, amsthm}
\usepackage[ruled,linesnumbered]{algorithm2e}
\usepackage{graphicx}
\usepackage{color}
\usepackage[top=3cm, bottom=2.5cm, left=2cm, right=2cm]{geometry}
\usepackage{float, caption, subcaption}

\input{epsf}

\newcommand{\bd}{\begin{description}}
\newcommand{\ed}{\end{description}}
\newcommand{\bi}{\begin{itemize}}
\newcommand{\ei}{\end{itemize}}
\newcommand{\be}{\begin{enumerate}}
\newcommand{\ee}{\end{enumerate}}
\newcommand{\beq}{\begin{equation}}
\newcommand{\eeq}{\end{equation}}
\newcommand{\beqs}{\begin{eqnarray*}}
\newcommand{\eeqs}{\end{eqnarray*}}

\definecolor{DarkGreen}{rgb}{0.2, 0.6, 0.3}

\newtheorem{theorem}{Theorem}[section]

\newtheorem{corollary}[theorem]{Corollary}
\newtheorem{case}{Case}

\newtheorem{claim}{Claim}

\newtheorem{proposition}{Proposition}[section]

\newtheorem{problem}{Problem}

\newtheorem{observation}{Observation}[section]
\begin{document}
\title{\textbf{Ramsey Achievement Games on Graphs : Algorithms and Bounds} \footnote{Supported by the National
Science Foundation of China (No. 12061059) and the Qinghai Key Laboratory of
Internet of Things Project (2017-ZJ-Y21).} }

\author{Xiumin Wang
\footnote{School of Mathematics and Statistics, Qinghai Normal
University, Xining, Qinghai 810008, China. {\tt wangxiumin1111@163.com}},
\  Zhong Huang\footnote{Corresponding author: School of Information and Mathematics, Yangtze
University, Jingzhou 434023, China. {\tt hz@yangtzeu.edu.cn}},
\  Xiangqian Zhou\footnote{Department of Mathematics and Statistics, Wright State University, Dayton, OH, 45435, USA. {\tt xiangqian.zhou@wright.edu}},
\  Ralf Klasing\footnote{Universit\'{e} de Bordeaux, Bordeaux INP, CNRS, LaBRI, UMR 5800, Talence, France.
{\tt ralf.klasing@labri.fr}},
\ Yaping Mao\footnote{Faculty of Environment and Information
Sciences, Yokohama National University, 79-2 Tokiwadai, Hodogaya-ku,
Yokohama 240-8501, Japan. {\tt maoyaping@ymail.com}}}
\date{}
\maketitle

\begin{abstract}
In 1982, Harary introduced the concept of Ramsey achievement game on graphs.
Given a graph $F$ with no isolated vertices. Consider the following game played on the complete
graph $K_n$ by two players Alice and Bob.
First, Alice colors one of the edges of $K_n$ blue,
then Bob colors a different edge red, and so on.
The first player who can complete the formation of $F$
in his color is the winner. The minimum $n$ for
which Alice has a winning strategy is the achievement number of $F$, denoted by $a(F)$.
If we replace $K_n$ in the game by the completed bipartite graph $K_{n,n}$, we get the bipartite achievement number, denoted by $\operatorname{ba}(F)$. In his seminal paper, Harary proposed an open problem of determining bipartite achievement numbers for trees. In this paper, we correct $\operatorname{ba}(mK_2)=m+1$ to $m$ and disprove $\operatorname{ba}(K_{1,m})=2m-2$ from Erickson and Harary, and
extend their results on bipartite achievement numbers. We also find the exact values of achievement numbers for matchings, and the exact values or upper and lower bounds of bipartite achievement numbers on matchings, stars, and double stars. Our upper bounds are obtained by deriving efficient winning strategies for Alice. \\[2mm]
\end{abstract}

\section{Introduction}

Combinatorics and graph theory play an important role in studying game theory; see \cite{ADMS21, Beck08, BFIN2021, BGGN18, CINP20, CNST23, CMINPS18, Kearns07}.  One good example is the domination game: in a graph, a vertex is said to dominate itself and its neighbors. The domination game is a two-player game played on a finite graph $G$ by two players,  who take turns choosing a vertex from $G$ such
that whenever a vertex is chosen by either player, at least one additional vertex is dominated. There had been a lot of research done on the domination game: see \cite{BKD2010, BIK2021, BDK2022}. For avoidance games based on graph theory, we refer to \cite{HHT03, FH89, Slany01}.

Ramsey theory, introduced in 1930 \cite{Ramsey}, is an important and yet very difficult subject in combinatorics. Roughly speaking, Ramsey theory aims to find certain monochromatic structures
for any coloring of a huge object \cite{GRS90}.

Given two graphs $G$ and $H$, the \emph{Ramsey number}
$r(G,H)$ is defined as the minimum
number of vertices $n$ needed so that every red/blue-edge-coloring of
$K_{n}$ contains a red copy of $G$ or a blue copy of $H$. If $G=H$, we write $r(G)=r(G,H)$ for short. Furthermore, if $G=K_n$, we write $r(n)=r(K_n)$. Given graphs $F,G,H$, does every edge-coloring of $F$ with red and blue contains either a red $G$ or a blue $H$? This problem is $coNP$ for fixed $G$ and $H$, $\prod ^p_2$-complete \cite{Schaefer99} in general.
The computational complexity of determining Ramsey numbers is still open. Some discussions about this problem can be found in \cite{Haanpaa00, Rosta04, Schweitzer09}. Algorithms for computing Ramsey numbers can be found in \cite{Overberghe20}.

Sah \cite{Sah} obtained the following result.
\begin{theorem}{\upshape \cite{Sah}}\label{th-Sah}
There is an absolute constant $c>0$ such that for $k\geq 3$,
$$
r(k+1)\leq {2k\choose k}e^{-c(\log k)^2}.
$$
\end{theorem}

Ramsey number has wide applications in the fields of communications,
information retrieval in computer science, and decision-making; see
\cite{Roberts, Rosta04} for examples. Interested reader may refer to \cite{MR1670625} for a dynamic survey and some papers \cite{GO22, Overberghe20} of small Ramsey numbers.

Motivated by Ramsey theory, Harray~\cite{HF1982} studied games related to Ramsey theory, which has application on economics \cite{SSV03}. Conversely, neural networks can be used in Ramsey theory \cite{GLZ22}. These games received more attention recently; see \cite{CDLM19,Gebauer,HNP17,MS10}.
In this paper we will study one such game called the achievement game. More specifically, we study the achievement game on a complete graph $K_n$ and the achievement game on the complete bipartite graph, introduced by Harary in \cite{HF1982} and by Erickson and Harary in \cite{HF1983}, respectively.

Given a finite graph $F$ with no isolated vertices,  the achievement game of $F$ on the complete graph $K_n$, denoted by $(F, K_n, +)$, is played by two players called Alice and Bob. In each round,  Alice first chooses an uncolored edge $e$ and color it blue, then Bob choose another uncolored edge $f$ and color it red;  the player who can first complete the formation of $F$ in his color is the winner.  The \emph{achievement number} of $F$, denoted by $a(F)$ is defined to be the smallest $n$ for which Alice has a winning strategy. Harary \cite{HF1982} presented the following results.

~\\

\begin{theorem}{\upshape
\cite{HF1982}}\label{Thm1} Let $F$ be a graph. Then
\begin{center}
\begin{tabular}{ccccccccccc}
\cline{1-11}
$F$ & $K_2$ & $P_3$ & $2K_2$ & $P_4$ & $K_{1,3}$ & $K_3$ & $C_4$&$K_3\cdot K_2$& $K_4-e$ & $K_4$\\[1mm]
\cline{1-11}
$a(F)$ & $2$ & $3$ & $5$ & $5$ & $5$ & $5$ & $6$&$5$& $7$& $10$\\[0.1cm]
$r(F)$ & $2$ & $3$ & $5$ & $5$ & $6$ & $6$ & $6$&$7$& $10$& $18$\\[0.1cm]
\cline{1-11}
\end{tabular}
\end{center}
\end{theorem}

For two bipartite graphs $G$ and $H$, the \emph{bipartite Ramsey number} ${\rm br}(G,H)$ of $G$ and $H$ is the smallest positive
integer $r$ such that every red-blue coloring of the $r$-regular complete bipartite graph $K_{r,r}$ results in either a red $G$ or a blue $H$. If $G=H$, then we write ${\rm br}(G)$. For more details on bipartite Ramsey number, we refer to \cite{GRS90}.

Harary \cite{HF1982} introduced the concept of bipartite achievement number.
Let $M$ be a bipartite graph with no isolated vertices. The achievement game of $M$ on the complete bipartite graph $K_{n,n}$, denoted by $(M, K_{n,n}, +)$, is played in a similar manner.  The \emph{bipartite achievement number} $\operatorname{ba}(M)$ is the minimum $n$ for which Alice has a winning strategy in the game $(M, K_{n,n}, +)$.

For small graphs, Erickson and Harary \cite{HF1983} obtained the following
results.
\begin{theorem}[\cite{HF1983}]\label{Thm2}
Let $M$ be a bipartite graph. Then
\begin{center}
\begin{tabular}{ccccccccccccc}
\cline{1-13}
$M$ & $K_2$ & $P_3$ & $2K_2$ & $P_4$ & $K_{1,3}$ & $P_3\cup K_2$ & $C_4$ & $P_5$ & $B_{3,2}$ & $K_{1,4}$ & $K_{2,3}-e$ & $K_{2,3}$ \\[1mm]
\cline{1-13}
$\operatorname{ba}(M)$ & $1$ & $2$ & $3$ & $3$ & $4$ & $3$& $4$ & $4$ & $5$ & $6$ & $4$ & $6$ \\[0.1cm]
$\operatorname{br}(M)$ & $1$ & $3$ & $3$ & $3$ & $5$ & $3$ & $5$ & $4$ & $5$ & $7$ & $5$ & $9$\\[0.1cm]
\cline{1-13}
\end{tabular}
\end{center}
\end{theorem}

The bipartite achievement numbers for stars, matchings and paths were also given in \cite{HF1983}.
\begin{theorem}[\cite{HF1983}]\label{ThmMH}
$(1)$ For $m\geq 2$,
 $$\operatorname{ba}(K_{1,m})=2m-2.$$
$(2)$ For $m\geq 2$,
$$
\operatorname{ba}(P_m)=
\begin{cases}
m-1, & m=2,3;\\
\lfloor\frac{m+3}{2}\rfloor, & m\geq 4.
\end{cases}
$$
$(3)$ For $m\geq 1$,
$$
\operatorname{ba}(mK_2)=
\begin{cases}
1, & m=1;\\
m+1, & m\geq 2.
\end{cases}
$$
\end{theorem}

Slany \cite{Slany01} considered combinatorial achievement games based on graph Ramsey theory: The players take turns in coloring edges of a graph $G$, each player being assigned a distinct color and choosing one so far uncolored edge per move. In achievement games, the first player that completes a monochromatic subgraph isomorphic to $A$ wins. He proved that general graph Ramsey achievement games are $PSPACE$-complete; see \cite{Slany01} for more details.

We present our main results in Section~\ref{sec:preli} after we introduce more notions and definitions. The proofs of our results are presented in Sections ~\ref{sec:complete} and \ref{sec:bipartite}.

\section{Preliminaries and Main results} \label{sec:preli}

Let $G$ be a graph. We will use $V(G)$, $E(G)$ to denote the vertex set and the edge set of $G$, respectively.
Furthermore, we will use $d_G(u)$ to denote the number of edges incident to vertex $u$ in graph $G$.
Let $E' \subset E(G)$. We will use $G[E']$ to denote the subgraph generated by the edge set $E'$. We will use $V[E']$ to denote the set of vertices incident with edges in the edge set $E'$.

The \emph{double star} $S_{m, n}$, where $n \geq m \geq 1$, is the graph consisting of the disjoint union of two stars $K_{1,m}$ and $K_{1,n}$ together with an edge joining their centers. Let $P_n$ denote the path on $n$ vertices.

Let $F$ be a graph (resp. bipartite graph)  and let $G$ be a complete graph (resp. complete bipartite graph), we will use $(F, G, +)$ to denote the achievement (resp. bipartite achievement)  game of $F$ played on $G$.  We will use $e_i$  (resp. $f_i$) to denote the edge chosen by Alice (resp. Bob) in the $i$-th round.

The following observations are immediate.
\begin{observation}\label{obs1}
Let $F$ be a graph. If $H$ is a spanning subgraph of $F$, then
$a(H)\leq a(F)$.
\end{observation}

\begin{observation}\label{obs2}
Let $F$ be a graph. Then $a(F)\leq r(F)$.
\end{observation}

\begin{corollary}\label{cor-a-lower}
Let $F$ be a graph of order $n$. Then
$$
n\leq a(F)\leq {2n-2\choose n-1}e^{-c(\log (n-1))^2}.
$$
Moreover, the lower bound is sharp.
\end{corollary}
\begin{proof}
It is clear that $a(F)\geq n$. From Observation~\ref{obs1}, we have $a(F)\leq a(K_n)$. From Theorem~\ref{th-Sah}, we have
$$
a(F)\leq a(K_n)\leq r(n)\leq {2n-2\choose n-1}e^{-c(\log (n-1))^2},
$$
where $n\geq 4$ and $c$ is an absolute constant.

The sharpness of the lower bounds follows immediately from Theorem~\ref{thm:mk2-1}, where we will show that $a(mK_2) =2m$ for $m \geq 3$.
\end{proof}

The following result shows that the graphs with large maximum degree can not reach the lower bound of Corollary~\ref{cor-a-lower}.
\begin{proposition}\label{cor-a-lower2}
Let $F$ be a graph of order $n$, if $\Delta(F)\geq n-2$, then $a(F)\geq n+1$.
\end{proposition}

The next observation is also straightforward.

\begin{observation}\label{obs3}
$(1)$ If $F$ is a bipartite graph, then $2\operatorname{ba}(F)\leq\operatorname{br}(F)$.

$(2)$ Let $F$ be a bipartite graph. If $H$ is a spanning subgraph of $F$, then
$\operatorname{ba}(H)\leq \operatorname{ba}(F)$.
\end{observation}

\begin{corollary}\label{cor-ba-lower}
If $F=(A,B)$ is a bipartite graph with two parts of order $|A|=a$ and $|B|=b$, then
$$
\max\{a,b\}\leq \operatorname{ba}(F)\leq \operatorname{br}(F)\leq \operatorname{br}(K_{a,b}).
$$
\end{corollary}

\begin{corollary}\label{cor-a-D}
Let $F=(A,B)$ be a bipartite graph with two parts of order $|A|=a$ and $|B|=b$. If $\operatorname{ba}(F)\geq \max\{a,b\}+1$, then $\Delta(F)\geq \max\{a,b\}-1$.
\end{corollary}

Harary \cite{HF1982} proposed the following problem.
\begin{problem}{\upshape \cite{HF1982}}\label{prob}
Determine bipartite achievement numbers for trees (especially stars and paths), for even cycles, and
for other families of bipartite graphs.
\end{problem}

Motivated by Problem \ref{prob}, we study the achievement numbers and bipartite achievement numbers of some special graphs.  Our main results include the following:

\begin{itemize}

\item[1.] We show that $a(mK_2) = 2m$ for $m \geq 3$ (Theorem~\ref{thm:mk2-1}) and $a(K_{1,n})\geq n+2$ for $n\geq 3$ (Theorem~\ref{thm:k1n-lower});

\item[2.]  Theorem~\ref{ThmMH} (Erickson and Harary~\cite{HF1983}) claims that $\operatorname{ba}(mK_2) = m+1$ for $m\geq 2$. We prove that this value is not accurate in Theorem~\ref{Thm3}, where we obtain that $\operatorname{ba}(mK_2)=m$ for $m\geq 4$.

\item[3.] We obtain the exact value of $\operatorname{ba}(S_{n,n})$ and an upper bound and lower bound for $\operatorname{ba}(S_{m,n})$ when $m\neq n$ (Theorem~\ref{Thm4});

\item[4.]  Theorem~\ref{ThmMH} (Erickson and Harary~\cite{HF1983}) claims that $\operatorname{ba}(K_{1,m}) = 2m-2$. We prove that this value is not accurate in Theorem~\ref{thm:ba(m-star)}, where we obtain that $m+1\leq \operatorname{ba}(K_{1,m})\leq 2m-3$ for $m\geq 3$.

\item[5.] Let $K_{1,m} \cup nK_2$ be the disjoint union of an $m$-star and an $n$-matching. We study $\operatorname{ba}(K_{1,m} \cup nK_2)$ in Theorem~\ref{the12} and find either the exact value or upper/lower bounds.

\item[6.] We study the bipartite achievement games $(F, K_{n_1, n_2}, +)$ where $n_1 \neq n_2$ and $F$ is isomorphic to $mK_2$, $S_{m,n}$, or $K_{1,m}$. We give sufficient conditions (some are necessary too) on when Alice has a winning strategy (Theorem~\ref{Thm8}, and Theorem~\ref{Thm7}).

\end{itemize}

\section{Results for the achievement numbers} \label{sec:complete}

In this section we will study the achievement numbers. Our first result is on $a(mK_2)$. We give the exact value of $a(mK_2)$ for all positive integer $m$.

\begin{theorem} \label{thm:mk2-1}
$$
\operatorname{a}(mK_2)=
\begin{cases}
5, & m=2;\\
2m, & m\neq 2.
\end{cases}
$$
\end{theorem}
\begin{proof}
Clearly the theorem holds for $m=1,2$ by Theorem \ref{Thm1}. We can assume that $m\geq 3$ in the following. Note that $|V(mK_2)| = 2m$. It follows that $a(mK_2)\geq 2m$.  To show that $a(mK_2)\leq 2m$, we will construct a winning strategy for Alice in the achievement game $(mK_2,K_{2m},+)$.

Next we will proof that Algorithm 1 implies an Alice's winning strategy.

\begin{algorithm}[H]\footnotesize
\caption{Alice's winning strategy in $(mK_2, K_{2m}, +)$}
\label{Ag1}

\KwIn {The achievement game $(mK_2, K_{2m}, +)$; Bob choose $f_i$ in the $i$-th round}
\KwOut {Alice choose $e_i$ in the $i$-th round}
$B\leftarrow \{uv\}$ \tcp{ $uv$ is an arbitrary edge of the graph $G=K_{2m}$}
$R\leftarrow\emptyset$\;
$e_1 \leftarrow uv$ \;
  \For{$i\leftarrow 2$ \KwTo $m-2$}
    {
       $e_{i} \leftarrow uv$ \tcp{ $u,v\in V(G)\setminus V(B)$ such that $d_{G[R]}(u)$ and $d_{G[R]}(v)$ are as large as possible}
       $B\leftarrow B\cup \{e_{i}\}$\;
       $R\leftarrow R\cup \{f_{i-1}\}$
    }
$e_{m-1} \leftarrow uv $\tcp{ $u,v\in V(G)\setminus V(B)$ such that $d_{G[R]}(u)$ and $d_{G[R]}(v)$ are as small as possible}
$wx \leftarrow f_{m-1}$\;
\eIf{$w,x\in V(G) \setminus V(B)$}
  {
    $e_{m} \leftarrow uv$ \tcp{ $u\in V(f_{m-1}), v\in V(e_{m-1})$ such that $d_{G[R]}(u)$ and $d_{G[R]}(v)$ are as large as possible}
    $B\leftarrow B\cup \{e_{m}\} \setminus \{e_{m-1} \}$\;
    $wx \leftarrow f_{m}$\;
    \eIf{ $w,x\in V(G) \setminus V(B)$ }
      {
        $e_{m+1} \leftarrow uv$ \tcp{ $u\in V(f_{m-1})\setminus V(e_{m})$ and $v\in V(e_{m-2})$}
        $B\leftarrow B\cup \{e_{m+1}\}\setminus \{e_{m-2}\}$\;
        $g_1\leftarrow uv$ \tcp{ $u,v \in V(G) \setminus V(B)$.}
        $B\leftarrow B\cup \{e_{m-1}\}\setminus \{e_{m}\}$\;
        $g_2\leftarrow uv$ \tcp{ $u,v \in V(G) \setminus V(B)$.}
        $e_{m+2} \leftarrow uv$ \tcp{ $uv \in  \{g_1,g_2\}\setminus \{f_{m+1}\}$}
      }
      {
        $e_{m+1} \leftarrow uv$ \tcp{ $u,v \in V(G) \setminus V(B)$}
      }
  }
  {
    $e_{m} \leftarrow uv$ \tcp{ $u,v \in V(G) \setminus V(B)$}
  }
\end{algorithm}

Recall that $e_i$ (resp. $f_i$) is the edge chosen by Alice (resp.  Bob) in the $i$-th round.
Let $B_i = \{e_1, e_2, \cdots, e_i\}, R_i = \{f_1, f_2, \cdots, f_i\}$ ($R_0=\emptyset$) and let $H_i$ be the subgraph spanned by $V(G) \setminus V(B_i)$.
For $1\leq i \leq m-2$, Alice chooses $e_{i} \in E(H_i)$ such that the degree of two endpoints of $e_{i}$ in $G[R_{i-1}]$ are as large as possible.
Then we have that $V(R_{m-3})\subset V(B_{m-2})$.
Alice choose $e_{m-1}$ such that the degree of two endpoints of $e_{m-1}$ in $G[R_{m-2}]$ are as small as possible.
If $f_{m-1} \notin E(H_{m-1})$, then Alice choose $e_{m}\in E(H_{m-1})$ and $B_m$ is a blue copy of $mK_2$.

We can assume that $f_{m-1} \in E(H_{m-1})$.
Alice choose $e_{m}$ such that one endpoint is in $V(f_{m-1})$ and the other endpoint is in $V(e_{m-1})$.
The degree of two endpoints of $e_{m}$ in $G[R_{m}]$ are as large as possible.
Let $M_0=B_m\setminus \{e_{m-1}\}$.
Then $R_m$ is not a matching by the choice of $e_{m-1}$.
If $V(f_{m})\neq V(G)\setminus V(M_0)$, then Alice choose $V(e_{m+1})= V(G)\setminus V(M_0)$ and $ M_0 \cup \{e_{m+1}\}$ is a blue copy of $mK_2$.

Next we can assume that $V(f_{m}) = V(G)\setminus V(M_0)$.
Alice choose $e_{m+1}$ such that one endpoint is in $V(f_{m-1})\setminus V(e_{m})$ the other endpoints is in $V(e_{m-2})$.
Clearly, $f_{m}$ is adjacent to $f_{m-1}$.
Then we can check that $R_{m+1}$ can not contain a red copy of $mK_2$.
Let $M_1=B_{m+1}\setminus \{e_{m-2},e_{m}\}, M_2=B_{m+1}\setminus \{e_{m-2},e_{m-1}\}, V(g_1)= V(G)\setminus V(M_1), V(g_2)= V(G)\setminus V(M_2)$.
Then $g_1, g_2\notin B_{m+1}\cup R_m$ by the choice of $\{e_{m-2}, e_{m-1}, e_m, e_{m+1}\}$.
Alice choose $e_{m+2}\in \{g_1,g_2 \} \setminus \{f_{m+1}\}$ and $B_{m+2}$ contains a blue copy of $mK_2$.
\end{proof}

Our second result in this section is the following lower bound for $a(K_{1,n})$.

\begin{algorithm}[H]\footnotesize
\caption{Alice's winning strategy in $(K_{1,n}, K_{n+1}, +)$}
\label{Ag2}
\KwIn {The achievement game $(K_{1,n}, K_{n+1}, +)$; Alice choose $e_i$ in the $i$-th round}
\KwOut {Bob choose $f_i$ in the $i$-th round}
$B\leftarrow \{e_1\}$ \;
$R\leftarrow\emptyset$\;
  \For{$i\leftarrow 1$ \KwTo $\lfloor\frac{E(G)}{2}\rfloor$}
    {
       \If{ $|V(G)\setminus V(R)|\geq 3$ }
         {
           $f_{i} \leftarrow wx$ \tcp{ $w,x\in V(G)\setminus V(R)$, $wx\neq e_i$ such that $d_{G[B]}(w)$ and $d_{G[B]}(x)$ are as large as possible}
         }
       \If{ $1\leq |V(G)\setminus V(R)|\leq 2$ }
         {
           $f_{i} \leftarrow wx$ \tcp{ $wx\in E(G)\setminus (B\cup R)$ and $w\in V(G)\setminus V(R)$}
         }
       \If{$|V(G)\setminus V(R)|=0$}
         {
           $f_{i} \leftarrow wx$ \tcp{ $wx\in E(G)\setminus (B\cup R) $}
         }
           $R\leftarrow R\cup \{f_{i}\}$\;
           $B\leftarrow B\cup \{e_{i+1}\}$
    }
\end{algorithm}

\begin{theorem}\label{thm:k1n-lower}
$a(K_{1,n})\geq n+2$ for $n\geq 3$.
\end{theorem}
\begin{proof}
Suppose that $V(K_{n+1}) = \{ v_1, v_2, \cdots, v_{n+1}\}$.  Recall that $e_i$ (resp. $f_i$) is the edge chosen by Alice (resp. Bob) in the $i$-th round, and $B_i$ (resp. $R_i$) is the set of blue edges (resp. red edges) after $i$ rounds. We will proof that Alice had no winning strategy in $(K_{1,n}, K_{n+1}, +)$. No matter how Alice choose edge $e_i$ in the $i$-th round, Algorithm 2 implies that there is no blue copy of $K_{1,n}$ when all the edges are colored.
We describe Algorithm 2 as follows.

In the $i$-th round for $1\leq i\leq  \lfloor\frac{E(G)}{2}\rfloor$.
\setcounter{case}{0}
\begin{case}
$|V(G)\setminus V(R_i)|\geq 3$.
\end{case}
In this case, Bob chooses $f_{i}$.
Clearly, $f_{i}$ can not be $e_i$.
The two endpoints of $f_{i}$ are in $V(G)\setminus V(R)$ such that $d_{G[B_i]}(w)$ and $d_{G[B_i]}(x)$ are as large as possible.
\begin{case}
$1\leq |V(G)\setminus V(R_i)|\leq 2$.
\end{case}
In this case, Bob chooses $f_{i}\in E(G)\setminus (B_{i+1}\cup R_i)$ such that one endpoint of $f_{i}$ is in $V(G)\setminus V(R)$.
\begin{case}
$|V(G)\setminus V(R_i)|=0$.
\end{case}
In this case, Bob chooses $f_{i}\in E(G)\setminus (B_{i+1}\cup R_i)$.

\begin{claim}\label{claimk1m}
In the $i$-th round, $|V(B_i)\setminus V(R_i)|\leq 1$ and the vertex $v_i\in V(B_i)\setminus V(R_i)$, if exists, is a leaf vertex in the graph $G[B_i]$.
\end{claim}
\begin{proof}[Proof of the Claim \ref{claimk1m}:]
We will proof by induction.
The claim holds for $i=1$.
Suppose that the claim holds for $i=k$.
\setcounter{case}{0}
\begin{case}
$|V(B_{k+1})\setminus V(R_k)|\geq 3$.
\end{case}
In this case, $|V(B_{k+1})\setminus V(R_k)|= 3$.
Every vertex in $V(B_{k+1})\setminus V(R_k)$ is a leaf vertex in graph $G[B_{i+1}]$.
Hence $|V(B_{k+1})\setminus V(R_{k+1})|=1$ and the claim holds for $i=k+1$.
\begin{case}
$|V(B_{k+1})\setminus V(R_k)|\leq 2$.
\end{case}
In this case, there exist at most one vertex in $V(B_{k+1})\setminus V(R_k)$ has degree at most 2 in graph $G[B_{i+1}]$.
Hence $|V(B_{k+1})\setminus V(R_{k+1})|\leq 1$ and the claim holds for $i=k+1$.
\end{proof}

If $u\in V(B_{i})\setminus V(R_i)$, then $d_{G[B_i]}(u)\leq 1$.
If $v\in V(B_{i+1})\setminus V(R_i)$, then $d_{G[B_{i+1}]}(u)\leq 2$.
Hence the edge $f_{i}$ can be chosen in each cases of Algorithm 2.
Now we have that $V(R_{\lfloor\frac{E(G)}{2}\rfloor})=V(G)$, and there is no blue copy of $K_{1,n}$.
\end{proof}

\section{Results for the bipartite achievement numbers} \label{sec:bipartite}

In this section we study achievement games played on a complete bipartite graph $K_{n_1, n_2}$. Throughout this section, we will use $X = \{x_i\,|\,1\leq i\leq n_1\}$ and $Y=\{y_i\,|\,1\leq i\leq n_2\}$ to denote the two parts of $K_{n_1,n_2}$. Recall that $e_i$ (resp. $f_i$) is the edge chosen by Alice (resp. Bob) in the $i$-th round.  We let $B_i =\{e_1, e_2, \cdots, e_i\}$,  $R_i = \{f_1, f_2, \cdots, f_i\}$,  $G_i = G[B_i]$ and $H_i = G[R_i]$.  We also use $F_i$ to denote the complete bipartite graph  $V(K_{n_1,n_2})\setminus V(G_i)$.

First we study the bipartite achievement numbers. And we give the exact value of $\operatorname{ba}(mK_2)$ for all positive integer $m$.
Note that, in the third part of Theorem~\ref{ThmMH}, Erickson and Harary \cite{HF1983} claim that $\operatorname{ba}(mK_2) =m+1$ for $m\geq 2$. However we show that in fact $\operatorname{ba}(mK_2)=m$ for $m\geq 4$ in Theorem~\ref{Thm3}.

\begin{theorem} \label{Thm3}
$$
\operatorname{ba}(mK_2)=
\begin{cases}
m+1, & \text{ $m=2$ or $m=3$};\\
m, & \text{ $m=1$ or $m\geq  4$} .
\end{cases}
$$
\end{theorem}

\begin{proof}
We can determine $\operatorname{ba}(mK_2)$ by exhaustion for $m\leq 4$.
Then we can assume that $m\geq 5$ in the following.
Since a monochromatic copy of $mK_2$ has at $2m$ vertices, then $\operatorname{ba}(mK_2)\geq m$.
To show that $\operatorname{ba}(mK_2)\leq m$, we will construct a winning strategy for Alice in the achievement game $(mK_2,K_{m,m},+)$.

Next we will proof that Algorithm 3 together with Algorithm 4 imply an Alice's winning strategy.
By similar argument in Theorem~\ref{thm:mk2-1}, we can check that $R_i$ can not be a red copy of $mK_2$ in the $i$-th round.
Then at the end of Algorithm 3 or Algorithm 4, $B_i$ contains a blue copy of $mK_2$.
\end{proof}

\begin{algorithm}[H]\footnotesize
\caption{Alice's winning strategy in $(mK_2, K_{m,m}, +)$ part one}
\label{Ag3}
\KwIn {The achievement game $(mK_2, K_{m,m}, +)$ for $m\geq 5$; Bob choose $f_i$ in the $i$-th round; $V(f_{m-2})\subseteq V(G)\setminus V(B_{m-2})$}
\KwOut {Alice choose $e_i$ in the $i$-th round}
$B\leftarrow \{uv\}$ \tcp{ $uv\in E(G)$ }
$R\leftarrow\emptyset$\;
$e_1 \leftarrow uv$ \;
  \For{$i\leftarrow 2$ \KwTo $m-3$}
    {
       $e_{i} \leftarrow uv$ \tcp{ $uv\in E(G)$, $v\in V(G)\setminus V(B)$ such that $d_{G[R]}(u)=0$ and $d_{G[R]}(v)$ are as large as possible}
       $B\leftarrow B\cup \{e_{i}\}$\;
       $R\leftarrow R\cup \{f_{i-1}\}$
    }
$e_{m-1} \leftarrow uv$ \tcp{ $uv\in E(G[V(G)\setminus B])$ such that $d_{G[R]}(u)=0$ and $v\in V(f_{m-2})$}
$B\leftarrow B\cup \{e_{m-1}\}$\;
\eIf{ $V(f_{m-1})= V(G)\setminus B$ }
      {
        $ab \leftarrow e_\alpha, cd \leftarrow e_\beta$ \tcp{ $e_\alpha,e_\beta \in B\setminus \{e_{m-1}\}$ such that $|E(G[V(\{e_\alpha, e_\beta\})])\cap R|\leq 1$ }
        $wx \leftarrow f_{m-1}$ \tcp{$a,c,u,w$ are in a same part of $G$ and $ad\notin R$}
        $e_{m} \leftarrow bw$\;
        \eIf{$f_{m}=ax$}
        {
          $e_{m+1} \leftarrow cx$\;
          $e_{m+2} \leftarrow yz$ \tcp{$yz\in\{ad,wd\}\setminus \{f_{m+1}\}$}
        }
        {$e_{m+1} \leftarrow ax$}
      }
      {
       $e_{m} \leftarrow uv$  \tcp{ $\{u,v\}=V(G)\setminus B$}
      }
\end{algorithm}

\begin{algorithm}[H]\footnotesize
\caption{Alice's winning strategy in $(mK_2, K_{m,m}, +)$ part two}
\label{Ag4}
\KwIn {The achievement game $(mK_2, K_{m,m}, +)$ for $m\geq 5$; Bob choose $f_i$ in the $i$-th round; $V(f_{m-2})\nsubseteq V(G)\setminus V(B_{m-2})$}
\KwOut {Alice choose $e_i$ in the $i$-th round}
$B\leftarrow \{uv\}$ \tcp{ $uv\in E(G)$ }
$R\leftarrow\emptyset$\;
$e_1 \leftarrow uv$ \;
  \For{$i\leftarrow 2$ \KwTo $m-3$}
    {
       $e_{i} \leftarrow uv$ \tcp{ $uv\in E(G)$, $v\in V(G)\setminus V(B)$ such that $d_{G[R]}(u)$ and $d_{G[R]}(v)$ are as large as possible}
       $B\leftarrow B\cup \{e_{i}\}$ \;
        $R\leftarrow R\cup \{f_{i-1}\}$
    }
$e_{m-1} \leftarrow uv$ \tcp{ $uv\in E(G[V(G)\setminus B])$ such that $d_{G[R]}(u)=d_{G[R]}(u)=0$}
$B\leftarrow B\cup \{e_{m-1}\}$\;
  \eIf{ $V(f_{m-1})= V(G)\setminus B$ }
      {
        $ab \leftarrow e_\alpha$\tcp{There exists $e_\alpha \in B\setminus \{e_{m-1}\}$ such that $E(G[V(\{e_\alpha, f_{m-1}\})])\cap R=\emptyset$ }
        $ uv\leftarrow e_{m-1}, wx\leftarrow f_{m-1}$ \tcp{ $a,u,w$ are in the same part of $G$ }
        $e_{m} \leftarrow bw$\;
          \eIf{$f_{m}=ax$}
           {
            $e_{m+1} \leftarrow ux$ \;
            $e_{m+2} \leftarrow yz$ \tcp{$yz\in\{av,aw\}\setminus \{f_{m+1}\}$}
           }
           {$e_{m+1} \leftarrow ax$}
     }
     {
       $e_{m} \leftarrow uv$  \tcp{ $\{u,v\}=V(G)\setminus B$ }
     }
\end{algorithm}

Recall that the double-star graph $S_{m,n}$ is obtained from $S_m$ and $S_n$ by adding an edge between the centers. Next we study the bipartite achievement number of $S_{m,n}$.

 \begin{theorem}\label{Thm4}
Let  $n \geq m \geq 1$. Then

$(1)$ if $m=n$, then $\operatorname{ba}(S_{m,n})= \operatorname{ba}(S_{n,n} )= 2n+1$; and

$(2)$ if $m\neq n$, then $n+1\leq \operatorname{ba}(S_{m,n})\leq 2n$.
\end{theorem}

\begin{proof}

First suppose that $m=n$.  To show $\operatorname{ba}(S_{n,n})\geq 2n+1$, we prove that Bob can stop Alice from forming $S_{n,n}$ on $K_{2n,2n}$.
If $e_i = xy$, then Bob would choose $f_i = x y^{\prime} $ where $y \neq y^{\prime}$. This is possible since every vertex in $X$ has even degree. It follows that every vertex in $X$ has blue degree at most $n$, and hence, Alice cannot form an $S_{n,n}$.

Now we show that Alice has a winning strategy on the game $(S_{n,n}, K_{2n+1, 2n+1}, +)$. Note that, Alice can form an $S_{n+1}$ in the (n+1)-th round regardless of Bob's moves. We may assume that $e_i = x_1y_i$ for $1 \leq i \leq n+1$. Note that, after $(n+1)$ rounds, there are exactly $n+1$ red edges chosen by Bob, so there exists $i_0 \in [1, n+1] $ such that $y_{i_0}$ is incident to at most one red edge. Then Alice will keep choosing edges incident to $y_{i_0}$ starting from the $(n+2)$-th round. It is clear that Alice will form an $S_{n,n}$ in the $(2n+1)$-th round.

Next we assume that $m < n$. Since $\Delta(S_{m,n}) = n+1$, it follows that $\operatorname{ba}(S_{m,n})\geq n+1$. To show that $\operatorname{ba}(S_{m,n}) \leq 2n$, it suffices to prove that Alice has a winning strategy in the game $(S_{m.n}, K_{2n, 2n}, +)$.
Assume that $e_1= x_1y_1$. Then either $x_1$ or $y_1$ is not incident to $f_1$. By symmetry, assume that $x_1$ is not incident to $f_1$. Now Alice will keep choosing edges incident to $x_1$ until the $(n+1)$-th round. So at the end of the $(n+1)$-th round, Alice forms a blue $S_{n+1}$ and there are exactly $n+1$ red edges (i.e, edges chosen by Bob). Without loss of generality, we may assume that $e_i = x_1y_i$ for $1 \leq i \leq n+1$.   Note that, at least one of $y_i$, $1 \leq i \leq n+1$ is incident to at most one red edge. By symmetry, assume that $y_1$ is incident to at most one red edge.  Then $y_1$ is incident to at least $2(n-1) \geq 2m$ uncolored edges. So starting from the $(n+2)$-th round, Alice will keep chosing edges incident to $y_1$. It is clear that Alice wins the game in the $(n+m+1)$-th round.  \end{proof}

Our next result studies the achievement game $(S_{m,n}, K_{n_1,n_2}, +)$ where $n_1$ may not be equal to $n_2$.

\begin{theorem}\label{Thm8}
Let $n \geq m \geq 1$ and consider the game  $(S_{m,n}, K_{n_1,n_2}, +)$.
\begin{itemize}
\item[1)] If $n=m$, then Alice has a winning strategy if and only if $\min\{n_1, n_2\} \geq 2n+1$;
\item[2)] If $m < n \leq 2m$ and Alice has a winning strategy, then one of $n_1$ and $n_2$ is at least $2m+1$ and the other one is at least $K$ for some $K \in [n+1, 2n]$;
\item[3)] If $n \geq 2m+1$, then Alice has a winning strategy if and only if $\min\{n_1,n_2\} \geq 2m+1$ and $\max\{n_1,n_2\} \geq 2n+1$.
\end{itemize}
\end{theorem}

\begin{proof}

For Part 1), first suppose that $n_2 = \min\{n_1, n_2\} \leq 2n$. Suppose that Alice chooses $e_i = x_iy_i$ in the ith round, then Bob will choose an uncolored edge incident to $x_i$, unless all edges incident to $x_i$ are colored. It follows that $d^r(x)  = d^b(x)$ or $d^b(x) -1$ for each $x \in X$.  Since $d^b(x) + d^r(x) \leq 2n$, every vertex in $X$ is incident to at most $n$ blue edges. Therefore, Alice can not form a $S_{n,n}$.

The 'if' part follows easily from Theorem~\ref{Thm4} 1).

For Part 2),  a similar argument shows that one of $n_1$ and $n_2$ is at least $2m+1$. So Part 2) follows easily from the following two claims:
\begin{claim}\label{claim1}
Alice has no winning strategy in $(S_{n,m}, K_{n,2m+1}, +)$
\end{claim}
\begin{proof}[Proof of the Claim~\ref{claim1}:]
Suppose that Alice chooses $e_i = x_iy_i$ in the ith round, then Bob will choose an uncolored edge incident to $x_i$, unless all edges incident to $x_i$ are colored.  It follows that every vertex in $X$ is incident to at most $m+1$ blue edges. Every vertex in $Y$ has degree $n$, and thus is incident to at most $n$ blue edges. So Alice can not form an $S_{n,m}$.
\end{proof}

\begin{claim}\label{claim2}
Alice has a winning strategy in $(S_{n,m}, K_{2n,2m+1}, +)$.
\end{claim}
\begin{proof}[Proof of the Claim~\ref{claim2}:]
Since every vertex $x \in X$ has degree $2m+1$, Alice can first form a $K_{1,m+1}$ centered $x_1$  after exactly $m+1$ rounds, regardless of Bob's move. By symmetry, we may assume that $e_i = x_1y_i$ for $1 \leq i \leq m+1$.   Note that, after $m+1$ rounds, there are exactly $m+1$ red edges chosen by Bob. If there exists $i_0 \in \{1, 2, \cdots, m+1\}$ such that $y_{i_0}$ is not incident any red edges, then Alice will keep choosing uncolored edges  incident to $y_{i_0}$ starting from the $(m+2)$-th round. Since $d(y_{i_0}=2n$, Alice can form at $S_{n+1}$ centered at $y_{i_0}$ in the $(m+n+1)$-th round, and hence, Alice forms an $S_{n.m}$.  So we may assume that after $m+1$ rounds, $y_i$ is incident to exactly one red edge, for each $1 \leq i \leq m+1$. Since Bob may stop Alice from winning by forming a red $S_{n,m}$, we may further assume that $f_i = x_2 y_i$ for $1 \leq i \leq m+1$.

Now Alice will choose another uncolored incident to $x_1$ in the $(m+2)$-th round, say $e_{m+2} = x_1y_{m+2}$.  If $f_{m+2}$ is incident to $x_1$, then clearly $f_{m+2}$ is not incident to $y_{m+2}$.  So Alice can keep choosing uncolored edges incident to $y_{m+2}$. Since  $d(y_{m+2})=2n$, Alice can form an $K_{1,n+1}$ centered at $y_{m+2}$ to win the game.  So $f_{m+2}$ is not incident to $x_1$. Now if $f_{m+2}$ is not incident to $y_{m+2}$, then Alice can keep choosing edges incident to $y_{m+2}$ to win the game. So Bob is forced to pick  $f_{m+2} = x_2y_{m+2}$. Alice then repeats the same strategy until Alice forms an $K_{1,n+1}$ centered at $x_1$. Now note that each $y_i$, $1 \leq i \leq n+1$, is incident to exactly one blue edge and one red edge. Since $d(y_1)= 2n \geq  2(m+1) > 2m+1$, Alice can keep choosing edges incident to $y_1$ to form an $S_{n,m}$.
\end{proof}

Part 3) follows easily from Claim~\ref{claim3} and Claim~\ref{claim4}.

\begin{claim}\label{claim3}
Alice has a winning strategy in $(S_{n,m}, K_{2n+1,2m+1}, +)$.
\end{claim}
\begin{proof} [Proof of the Claim~\ref{claim3}:]   Assume that $e_1 =x_1y_1$. Since $d (x_1)=2m+1$, Alice can keep choosing edges incident to $x_1$ in the first $m+1$ rounds, regardless of Bob's move.   So at the end of the $(m+1)$-th round, Alice forms a blue $S_{m+1}$ and there are exactly $m+1$ red edges, (edges chosen by Bob). Assume that $e_i = x_1 y_i$ for $1 \leq i \leq m+1$. Note that there exists $j$ with $1 \leq j \leq m+1$ and $y_j$ is incident to at most one red edges. Assume by symmetry that $j=1$. Now Alice will keep choosing uncolored edges incident to $y_1$ starting from the $(m+2)$-th round. It is clear that Alice will form an $S_{n,m}$ in the $(n+m+1)$-th round. \end{proof}

\begin{claim}\label{claim4}
Alice has no winning strategy in $(S_{n,m}, K_{2n,2m+1}, +)$.
\end{claim}

\begin{proof} [Proof of the Claim 4:]   Assume that $e_i = x_iy_i$. Then Bob will simply choose an uncolored edge incident to $y_i$. Since $d(y_i)=2n$, it follows that every vertex in $Y$ is incident to at most $n$ blue edges. So Alice can not form an $S_{n,m}$. \end{proof}\end{proof}

The next result gives the exact value of $\operatorname{ba}(K_{1,4})$.

\begin{theorem}\label{theK_{1,4}}
$\operatorname{ba}(K_{1,4})=5$.
\end{theorem}

\begin{proof}

First we show that Alice does not have a winning strategy in $(K_{1,4}, K_{4,4}, +)$.  Suppose that in the i-th round, Alice picks $e_i = xy$. If both $d^r(x)=0$ and $d^r(y)=0$, then we pick $u \in \{x, y\}$ such that $d^b(u) = \max \{d^b(x), d^b(y)\}$; without loss of generality, assume that $u=x$. Then we pick $v \in Y$ such that $uv$ is uncolored and $d^r(v)$ is as small as possible, let $f_i = uv$. If exactly one of $d^r(x)$ and $d^r(y)$ is 0, say $d^r(x)=0$, then pick $v \in Y$ such that  $xv$ is uncolored and $d^r(v)$ is as small as possible; and let $f_i = xv$. Finally if neither $d^r(x)=0$ nor $d^r(y)=0$, then Bob pick $f_i = x^{\prime}y^{\prime}$ such that $x^{\prime} \in X$, $y^{\prime}  \in Y$,  and each of $x^{\prime}$ and $y^{\prime}$ is incident to as many blue edges and as few red edges as possible. It is straightforward to check that the above procedure produces a strategy for Bob to stop Alice from winning in $(K_{1,4}, K_{4,4}, +)$.

Next we show that Alice has a winning strategy in $(K_{1,4}, K_{5,5}, +)$.  Assume that $e_1 = x_1y_1$. By symmetry, either $f_1= x_1y_2$ or $f_1 = x_2y_2$. First assume that  $f_1 = x_1y_2$. Then Alice choose $e_2 = x_2y_1$; this would force $f_2 = x_3y_1$. Now Alice chooses $e_3 = x_4y_1$, forcing $f_3 = x_5y_1$. Then Alice chooses $e_4 = x_2y_2$, which would force $f_4=x_2y_3$. So then Alice chooses $e_5= x_2y_4$ and forcing $f_5 = x_2y_5$. Next Alice chooses $e_6=x_4y_4$ and each of $x_4$ and $y_4$ is incident to two blue edges and no red edges.  It is clear that Alice can form a $K_{1,4}$ centered at $x_4$ or $y_4$, regardless of Bob's move.

If $f_1 = x_2y_2$, then Alice will choose $e_2= x_1y_2$, and hence this case is reduced to the previous case. \end{proof}

Note that, in the first part of Theorem~\ref{ThmMH}, Erickson and Harary \cite{HF1983} claim that $\operatorname{ba}(K_{1,m}) =2m-2$. We show that in fact the value of $\operatorname{ba}(K_{1,m})$ should be no more than $2m-3$. However, we could not decide the exact value.

\begin{theorem} \label{thm:ba(m-star)}
For $m\geq 3$, $m+1\leq \operatorname{ba}(K_{1,m})\leq 2m-3$.
\end{theorem}

\begin{proof}
First we show that Alice has a winning strategy in $(K_{1,m}, K_{2m-3, 2m-3}, +)$. Assume that $e_1 = x_1y_1$.  Then by symmetry, either $f_1=x_1y_2$ or $f_1= x_2y_2$. If $f_1 = x_1y_2$, then Alice chooses $e_2 = x_2y_1$. Now if $f_2$ is not incident to $y_1$, then Alice can keep chooses edges incident to $y_1$ to obtain a blue $K_{1,m}$. So $f_2 = x_3y_1$. Next  Alice chooses $e_3=x_4y_1$ and again forces $f_3=x_5y_1$. In the following step, Alice chooses $e_4 = x_2y_2$ and forces $f_4=x_2y_3$. Then Alice chooses $e_5=x_2y_4$ and forces $f_5=x_2y_5$. Now Alice chooses $e_6=x_4y_4$. Note that, each of $x_4$ and $y_4$ is incident to two blue edges and no red edges. Since $f_6$ is not incident to at least one of $x_4$ or $y_4$, say $x_4$. Then Alice will keep choosing edges incident to $x_4$ to form a $K_{1,m}$ first.  If Bob chooses $f_1=x_2y_2$, then Alice chooses $e_2=x_1y_2$, and hence it is reduced to the previous case by symmetry.

Next we show that Alice has no winning strategy in $(K_{1,m}, K_{m, m}, +)$.  Let $B_i$ (resp. $R_i$) be the set of blue (resp. red) edges after $i$-round.  It suffices to show that Bob can make his moves such that $|V(G[B_i]) \setminus V(G[R_i])| \leq 1$ and the vertex in $V(G[B_i]) \setminus V(G[R_i])$, if exists, is a leaf vertex of $G[B_i]$, for each $i\geq 1$.

We will proceed by induction. Clearly in the first round, Bob can choose $f_1$ such that $f_1$ is adjacent to $e_1$, so the statement holds for $i=1$.  Suppose that the statement holds for $i$, we will show that it holds for $i+1$. Clearly $|V(G[B_i])\setminus V(G[R_i])| \leq |V(G[B_{i+1}])\setminus V(G[R_i])|\leq |V(G[B_{i}])\setminus V(G[R_i])|+2$. So We have the following six cases:

\setcounter{case}{0}
\begin{case}
$|V(G[B_{i}])\setminus V(G[R_i])|=0$ and $|V(G[B_{i+1}])\setminus V(G[R_i])|=0$.
\end{case}
In this case, Bob chooses an arbitrary edge $f_{i+1}$ on the $(i+1)$th round.
Then $|V(G[B_{i+1}])\setminus V(G[R_{i+1})|=0$.

\begin{case}
$|V(G[B_{i}])\setminus V(G[R_i])|=0$ and $|V(G[B_{i+1}])\setminus V(G[R_i])|= 1$.
\end{case}
In this case, Bob chooses an arbitrary edge $f_{i+1}$ on the $(i+1)$th round.
Then $|V(G[B_{i+1}])\setminus V(G[R_{i+1}])|= 1$ and clearly the vertex in $V(G[B_{i+1}])\setminus V(G[R_{i+1}])$ is a leaf vertex in $G[B_i]$.

\begin{case}
$|V(G[B_{i}])\setminus V(G[R_i])|=0$ and $|V(G[B_{i+1}])\setminus V(G[R_i])|=2$.
\end{case}
In this case, suppose that $e_{i+1}=v_{1,i+1}v_{2,i+1}$. Then $v_{1,i+1},v_{2,i+1}\notin V(G[R_i])$.
So Bob chooses $f_{i+1}$ such that $f_{i+1}$ is incident with $v_{1,i+1}$.
Then $|V(G[B_{i+1}])\setminus V(G[R_{i+1}])|=1$ and $v_{2,i+1}$ is a leaf vertex in the graph $G[B_{i+1}]$.

\begin{case}
$|V(G[B_{i}])\setminus V(G[R_i])|=1$ and $|V(G[B_{i+1}])\setminus V(G[R_i])|= 1$.
\end{case}
In this case, suppose that $v_{0,i} \in G[B_{i}])\setminus V(G[R_i])$. Then $v_{0,i}$ is a leaf vertex in $G_i$.
There are at least three edges incident with $v_{0,i}$ in $K_{m,m}$.
Note that there are at most two edges incident with $v_{0,i}$ in $G[B_i]$. So Bob can choose $f_{i+1}$ such that $f_{i+1}$ is incident with $v_{0,i}$ in the $(i+1)$-th round.
Then $|V(G[B_{i+1}])\setminus V(G[R_{i+1}])|=0$.

\begin{case}
$|V(G[B_{i}])\setminus V(G[R_i])|=1$ and $|V(G[B_{i+1}])\setminus V(G[R_i])|=2$.
\end{case}
In this case, suppose that $e_{i+1}=v_{1,i+1}v_{2,i+1}$.
Then at least one of the endpoints of $e_{i+1}$ is not in $V(G_i)$,  say $v_{1,i+1}$.  Let $v_{0,i}$ be the vertex in $V(G[B_{i}])\setminus V(G[R_i])$. Then $v_{0,i}$ is a leaf vertex in the graph $G[B_i]$.
There are at least three edges incident with $v_{0,i}$ in $K_{m,m}$.
Note that there are at most two edges incident with $v_{0,i}$ in $G_i$. So Bob can choose $f_{i+1}$ such that $f_{i+1}$ is incident with $v_{0,i}$ in the $(i+1)$-th round.
Then $|V(G[B_{i+1}])\setminus V(G[R_{i+1}])|=1$ and $v_{1,i+1}$ is a leaf vertex in graph $G[B_{i+1}]$.

\begin{case}
$|V(G[B_{i}])\setminus V(G[R_i])|=1$ and $|V(G[B_{i+1}])\setminus V(G[R_i])|=3$.
\end{case}
In this case, suppose that $e_{i+1}=v_{1,i+1}v_{2,i+1}$.
Then $v_{1,i+1},v_{2,i+1}\notin V(G[B_i])$ and $V(G[B_i])\setminus V(G[R_i])=\{v_{0,i}\}$.
Bob chooses $f_{i+1}=v_{0,i}v_{1,i+1}$ in the $(i+1)$-th round.
Then $|V(G[B_{i+1}])\setminus V(G[R_{i+1}])|= 1$ and $v_{2,i+1}$ is a leaf vertex in the graph $G[B_{i+1}]$.

\end{proof}

Our next result studies the achievement game of $K_{1,m}$ played on $K_{n_1, n_2}$ where $n_1 \geq n_2$.

\begin{theorem}\label{Thm7}

Consider the game $(K_{1,m}, K_{n_1, n_2}, +)$ where $n_1 \geq n_2$.

\begin{itemize}
\item[1)] If $1 \leq n_2 \leq m-1$, then Alice has a winning strategy if and only if $n_1 \geq 2m-1$.

\item[2)] If $n_1 = 2m-2$, then Alice has a winning strategy if and only if $n_2 \geq m$.

\item[3)] If $n_1=n_2=2m-3$, then Alice has a winning strategy.

\end{itemize}
\end{theorem}

\begin{proof}

For part 1), first we show that Alice has no winning strategy on $(K_{1,m}, K_{2m-2, n_2}, +)$. Since $n_2 \leq m-1$, every $K_{1,m}$-subgraph must be centered in $Y$. Suppose that $e_i = xy$, then Bob would simply choose $f_i$ to be incident with $y$. This is possible since every $y \in Y$ has even degree (of $2m-2$). It is clear that Alice can not form a $K_{1,m}$.

Now suppose that $n_1 \geq 2m-1$ and  $e_1=x_1y_1$. Then it is clear that Alice can win if he keeps choosing uncolored edges incident to $y_1$.

For part 2), it is easy to see that Alice has no winning strategy on $(K_{1,m}, K_{2m-2, m-1}, +)$. We now show that Alice has a winning strategy on $(K_{1,m}, K_{2m-2, m}, +)$.    Assume that $e_1=x_1y_1$. Note that $f_1$ must be incident to $y_1$; as otherwise, Alice wins by repeatedly choosing uncolored edges incident to $y_1$. Assume that $f_1=x_2y_1$. Alice then chooses $e_2=x_3y_1$, again forcing $f_2=x_4y_1$. Alice will repeat the same strategy for $m-1$ rounds and forces Bob to choose edges incident to $y_1$. After all edges incident to $y_1$ are colored, Alice then chooses $e_m = x_1y_2$ and use the same strategy to force Bob choose edges incident to $y_2$ in the next $(m-1)$ round. After Alice repeat the same pattern for each of $y_1$, $y_2$, $y_{m-2}$. Then in the $(m-1)(m_2)+1$-th round, Alice chooses $x_1y_{m-1}$ and forms a blue $K_{1,m-1}$ centered at $x_1$. So Bob must choose $x_1y_m$ and then Alice can win by repeatedly choosing uncolored edges incident to $y_{m-1}$.

Part 3) follows from Theorem~\ref{thm:ba(m-star)}. \end{proof}

Next we study the graph $K_{1,m} \cup nK_2$, the disjoint union of a star and a matching, and prove the following result on the bipartite achievement numbers $\operatorname{ba}(K_{1,m} \cup nK_2)$ where $m \geq 1$ and $n \geq 1$.

\begin{theorem}\label{the12}

Suppose that $m, n \geq 1$. Then

\begin{itemize}
\item[1)] $\operatorname{ba}(K_{1,1} \cup nK_2) = \begin{cases}
n+2, & \text{ $n=1$ or $n=2$};\\
n+1, & \text{ $n=0$ or $n\geq 3$} .
\end{cases}
$;
\item[2)] if $2 \leq m \leq 4$ or $n \geq m-2$, then $\operatorname{ba}(K_{1,m} \cup nK_2)=m+n$;
\item[3)] if $m \geq 5$ and $n \leq m-3$, then $m+n \leq \operatorname{ba}(K_{1,m} \cup nK_2) \leq 2m-2$.

\end{itemize}
\end{theorem}
\begin{proof}

For part 1), note that $K_{1,1} \cup nK_2$ is isomorphic to  $(n+1) K_2$. Therefore, by Theorem~\ref{Thm3},  $
\operatorname{ba}(K_{1,1}\cup nK_2) = \operatorname{ba}((n+1)K_2)=
\begin{cases}
n+2, & \text{ $n=1$ or $n=2$};\\
n+1, & \text{ $n=0$ or $n\geq 3$} .
\end{cases}
$

Since $K_{1,m} \cup nK_2$ is a bipartite graph with one side having $m+n$ vertices, it follows that $\operatorname{ba}(K_{1,m}\cup nK_2)  \geq m+n$.

For Part 2), first assume that $m=2$. Since $K_{1,2} \cup nK_2$ is a bipartite graph with $n+1$ vertices on one side and $n+2$ vertices on the other side, it is clear that $\operatorname{ba}(K_{1,2}\cup nK_2)\geq n+2$. Now we show Alice has a winning strategy on $(K_{1,2} \cup nK_2, K_{n+2,n+2}, +)$. Note that  $\operatorname{ba}((n+1)K_2)=n+2$ and Alice can first form an $(n+1)K_2$ in the $(n+1)$-th round. At the end of the $(n+1)$-th round, there are exactly $n+1$ red edges and exactly two vertices, say $x\in X$ and $y\in Y$, not incident to any blue edges. So there are at least $2(n+1) - (n+1) = n+1$ uncolored edges and each of them can be chosen by Alice to form a $K_{1,2} \cup nK_2$ in the $(n+2)$-th round.

Next we show that Alice has a winning strategy in $(K_{1,4}\cup K_2, K_{5,5}, +)$.  Suppose that Alice chooses $e_1 =x_1y_1$. Then by symmetry either $f_1=x_1y_2$ or $f_1=x_2y_2$. In the former case, Alice chooses $e_2=x_2y_2$; while in the latter case, Alice chooses $e_2 = x_2y_1$. Now the two cases are symmetric to each other. So we may assume that  $f_1=x_1y_2$ and $e_2=x_2y_1$. Then Bob would choose $f_2=x_3y_1$ as otherwise Alice would form a $K_{1,4}$ in the fourth round and it is easy to check that Alice can form a $K_{1,4}\cup K_2$ in the fifth round.  Similarly, Alice chooses $e_3 = x_4y_1$ and forces $f_3= x_5y_1$.  Next Alice chooses $e_4=x_2y_2$, then we have the following four cases depending on $f_4$:
\setcounter{case}{0}
\begin{case}
$f_4=x_3y_3$.
\end{case}

Then Alice pick $e_5=x_2y_3$. Now Bob either pick an edge incident to $x_2$ to stop Alice from forming a $K_{1,4}$ centered at $x_2$, or pick an edge incident o $x_3$ and try to form a red $K_{1,4}$.  First assume that $f_5=x_2y_4$. Then Alice would choose $e_6= x_2y_5$. Note that each uncolored edge incident to $y_4$ can be chosen by Alice to form a $K_{1,4} \cup K_2$ and Bob can not stop it.  Next assume that $f_5=x_3y_4$. Then Alice must choose $e_7=x_3y_5$, forcing $e_7=x_2y_4$. Now Alice chooses $e_8=x_2y_5$. Note that, Bob can not win in the 8th round, but either $x_4y_4$ or $y_4y_5$ can be chosen by Alice in the 9th round to win the game.

\begin{case}
$f_4= x_3y_2$, $x_4y_4$, or $x_4y_2$.
\end{case}

In this case, Bob will need at least three more edges to form an $K_{1,4} \cup K_2$. Alice will choose two edges incident to $x_2$ in the next two round such that they are adjacent to as many red edges as possible, It is easy to check that Alice will win in the 7th round.

\begin{case}
$f_4=x_2y_3$.
\end{case}

In this case, Alice chooses $e_5=x_2y_4$, forcing $f_5= x_2y_5$.  Then Alice chooses $e_6 = x_4y_4$.  Note that, at this stage, Bob would need at least three more edges to form $K_{1,4} \cup K_2$.  Also assume by symmetry that $f_6 \neq x_4y_5$. So Alice chooses $e_7 = x_4y_5$. It is clear that Alice will win in the 8th round.

\begin{case}
$f_4=x_1y_3$.
\end{case}

In this case, Alice will pick $e_5 = x_2y_3$. Then if $f_5$ is not incident to $x_1$ or $x_2$, then Alice can win in two more moves while Bob cannot.  So we may assume that $f_5=x_1y_4$ or $f_5=x_2y_4$. If  $f_5=x_1y_4$, then $e_6=x_1y_5$. It is easy to check that Bob would need at least three more edges to win; while Alice only need two more edges to win regardless of Bob's choices. If $f_5=x_2y_4$, then Alice chooses $e_6=x_2y_5$. It is clear that Alice can win in the 7th round.

To complete the proof for part 2), it suffices to show that if $m \geq 3$ and $n \geq m-2$, then Alice has a winning strategy in $(K_{1,m} \cup nK_2, K_{m+n, m+n}, +)$.   Since $n \geq m-2$, $m+n \geq 2m-2$. Assume that $e_1=x_1y_1$ and that $f_1$ is not incident to $y_1$. By Symmetry, $f_1=x_1y_2$ or $f_1=x_2y_2$. In both cases, Alice chooses $e_2=x_2y_1$. Then either $f_2$ is adjacent to a blue edge or not. In the former case, Alice will pick another uncolored edge incident to $y_1$, in the latter case, assume that $f_2$ is incident to $x_3$, then Alice picks $e_3=x_3y_1$. Alice will repeat the same pattern for $m$-round, it is easy to see that at the end of the $m$-th round, Alice will form a $K_{1,m}$ and there is at most one red edge not adjacent to a blue edge.  Let $B$ be the set of vertices incident to a blue edge.  Then by deleting all vertices in $B$ and at most one other vertex in $Y$, we get a $K_{n, m+n-2}$ that contains no colored edges. Since $m+n-2 \geq n+1$ and Alice has a winning strategy in $(nK_2, K_{n.n+1}, +)$, Alice can first construct a $K_{1,m} \cup nK_2$ is $K_{n+m,n+m}$.

For part 3), we assume that $m \geq 5$ and $n \leq m-3$, we will show that Alice has a winning strategy on  $(K_{1,m}\cup nK_2$, $K_{2m-2,2m-2}, +)$. Assume that $e_1=x_1y_1$. Alice will use a similar strategy used in the previous paragraph to construct a $K_{1,m}$ in the $m$-th round centered at $y_1$; moreover, at the end of the $m$-th round, there is at most one red edge not adjacent  to any blue edges. By deleting all vertices incident to a blue edge and at most one other vertex in $Y$, we get a $K_{m-2, 2m-4}$ that contains no colored edges. Note that $\min\{m-2, 2m-4\} \geq n+1$. Since Alice has a winning strategy in $(nK_2, K_{n,n+1}, +)$, Alice can first construct a $K_{1,m} \cup nK_2$ is $K_{2m-2,2m-2}$.
\end{proof}

\section{Concluding remark}

In this paper, we find the exact values of achievement numbers for matchings, and the exact values or upper and lower bounds of bipartite achievement numbers on stars, matchings, and some of our proofs are obtained by deriving efficient algorithms. We study the Ramsey achievement games on graphs. It is interesting to study other graph classes and improve the bounds in Corollaries \ref{cor-a-lower} and \ref{cor-ba-lower}.


\begin{thebibliography}{1}

\bibitem{ADMS21}
E.C. Akrida, A. Deligkas, T. Melissourgos, P.G. Spirakis,
Connected subgraph defense games, \emph{Algorithmica} 83(11) (2021), 3403--3431.

\bibitem{Beck08}
J. Beck, \emph{Combinatorial games: tic-tac-toe theory}, Cambridge University Press, 2008.

\bibitem{BFIN2021}
J. Bensmail, F. Fioravantes, F.M. Inerney, N. Nisse,
The largest connected subgraph game, \emph{Algorithmica} 84(9) (2022), 2533--2555.

\bibitem{BGGN18}
B. Bosek, P. Gordinowicz, J. Grytczuk, N. Nisse, Joanna Sok\'{o}{\l}, M. \'{S}leszy\'{n}ska-Nowak, Localization game on geometric and planar graphs, \emph{Discrete Appl. Math.} 251 (2018), 30--39.


\bibitem{BKD2010}
B. Bre\u{s}ar, S. Klav\u{z}ar, D.F. Rall, Domination game and an imagination strategy, \emph{SIAM J. Discrete Math.} 24(3) (2010), 979--991.


\bibitem{BDK2022}
C. Bujt\'{a}s, P. Dokyeesun, S. Klav\u{z}ar, Thresholds for the monochromatic clique transversal game, \emph{Expo. Math.}, in press.


\bibitem{BIK2021}
C. Bujt\'{a}s, V. Ir\u{s}i\u{c}, S. Klav\u{z}ar, K. Xu, On Rall's $1/2$-conjecture on the domination game, \emph{Quaest. Math.} 44(12) (2021), 1711--1727.


\bibitem{CINP20}
N. Cohen, F.M. Inerney, N. Nisse, S. P\'{e}rennes,
Study of a combinatorial game in graphs through linear programming, \emph{Algorithmica} 82(2) (2020), 212--244.


\bibitem{CNST23}
M. Chudnovsky, S. Norin, P. Seymour, J. Turcotte, Cops and robbers on $P_5$-free graphs, \emph{arXiv preprint}
https://doi.org/10.48550/arXiv.2301.13175, 2023.


\bibitem{CMINPS18}
N. Cohen, N.A. Martins, F.M. Inerney, N. Nisse, S. P\'{e}rennes, R. Sampaio, Spy-game on graphs: Complexity and simple topologies, \emph{Theor. Comput. Sci.} 725 (2018), 1--15.


\bibitem{CDLM19}
D. Conlon, S. Das, J. Lee, T. M\'{e}sz\'{a}ros, Ramsey games near
the critical threshold, \emph{Random Struct. Algor.} 57(4) (2022), 940--957.


\bibitem{HF1983}
M. Erickson, F. Harary, \emph{Generalized Ramsey theory XV:
Achievement and avoidance games for bipartite graphs}, Graph Theory
Singapore 1983, Springer, Berlin, Heidelberg, 1984. 212--216.


\bibitem{FH89}
A.S. Fraenkel, F. Harary, Geodetic contraction games on graphs, \emph{Inter. J. Game Theory} 18 (1989), 327--338.


\bibitem{Gebauer}
H. Gebauer, Size Ramsey number of bounded degree graphs for games,
\emph{Combin. Probab. Comput.} 22(04) (2013), 499--516.

\bibitem{GLZ22}
A. Ghose, A. Levi, Y. Zhang, Graph neural networks for Ramsey graphs,
36th Conference on Neural Information Processing Systems (NeurIPS 2022) Workshop on MATH-AI.

\bibitem{GO22}
J. Goedgebeur, S.V. Overberghe, New bounds for Ramsey numbers $R(K_k-e,K_l-e)$, \emph{Discrete Appl. Math.}
307 (2022), 212--221.

\bibitem{GRS90}
R.L. Graham, B.L. Rothschild, J.H. Spencer, \emph{Ramsey Theory},
JOHN WILEY \& SONS, 1990.

\bibitem{Haanpaa00}
H. Haanp\"{a}\"{a}, \emph{Computational methods for Ramsey numbers}, Helsinki University of Technology, 2000.

\bibitem{HF1982}
F. Harary, \emph{Achievement and avoidance games for graphs},
North-Holland Mathematics Studies, Vol. 62, North-Holland, 1982.
111--119.

\bibitem{HHT03}
T.W. Haynes, M. Henning, C.A. Tiller, Geodetic achievement and avoidance games for graphs, \emph{Quaest. Math.}
26(4) (2003), 389--397.

\bibitem{HNP17}
D. Hefetz, C. Kusch, L. Narins, A. Pokrovskiy, C. Requil\'{e}, A.
Sarid, Strong Ramsey games: Drawing on an infinite board, \emph{J.
Combin. Theory, Ser. A} 150 (2017), 248--266.

\bibitem{Kearns07}
M. Kearns, Graphical Games. In Vazirani, Vijay V.; Nisan, Noam; Roughgarden, Tim; Tardos, \'{E}va (2007). Algorithmic Game Theory (PDF). Cambridge, UK: Cambridge University Press.



\bibitem{MS10}
M. Marciniszyn, R. Sp\"{o}hel, Online vertex-coloring games in
random graphs, \emph{Combinatorica} 30(1) (2010), 105--123.

\bibitem{Overberghe20}
S.V. Overberghe, \emph{Algorithms for Computing Ramsey Numbers}, Master of Science, Ghent University, 2020.

\bibitem{MR1670625}
S.~P. Radziszowski, Small {R}amsey numbers, \emph{Electron. J.
Combin.}, Dynamic Survey 1, 30 pp. (electronic), 1994.

\bibitem{Ramsey}
F.P. Ramsey, On a problem of formal logic, \emph{Proc. London. Math.
Soc.} (2)30 (1930), 264--286.

\bibitem{Roberts}
F.S. Roberts, Applications of Ramsey theory, \emph{Discrete Appl.
Math.} 9(3) (1984), 251--261.

\bibitem{Rosta04}
V. Rosta, Ramsey theory applications, \emph{Electron. J. Combin.}
11(1) (2004), 89.


\bibitem{Sah}
A. Sah, Diagonal Ramsey via effective quasirandomness, \emph{Duke Math. J.} 172(3) (2023), 545--567.

\bibitem{Schaefer99}
M. Schaefer, \emph{Graph Ramsey theory and the polynomial hierarchy}, Proceedings of the thirty-first annual ACM symposium on Theory of computing, ACM, 1999.

\bibitem{Schweitzer09}
P. Schweitzer, Problems of unknown complexity: graph isomorphism and Ramsey theoretic numbers, Diss. Saarbr\"{u}cken, Univ., Diss., 2009, 2009.

\bibitem{SSV03}
E. Shmaya, E. Solan, N. Vieille, An application of Ramsey theorem to stopping games, \emph{Games and Economic Behavior} 42(2) (2003), 300--306.

\bibitem{Slany01}
W. Slany, \emph{The Complexity of Graph Ramsey Games}. In: Marsland, T., Frank, I. (eds) Computers and Games. CG 2000. Lecture Notes in Computer Science, vol 2063. Springer, Berlin, Heidelberg, 2001.






\end{thebibliography}
\end{document}